\documentclass[12pt,twoside]{amsart}
\usepackage{amssymb}
\usepackage{amscd}

\title{Base point free theorem of Reid-Fukuda type}  
\author{Osamu Fujino} 
\date{\today}
\address{Research Institute for Mathematical Sciences\\ 
 Kyoto University, Kyoto 606-8502 Japan}
\email{fujino@kurims.kyoto-u.ac.jp}
\newcommand{\bQ}[0]{{\mathbb Q}}

\newcommand{\xDiff}[0]{{\operatorname{Diff}}}
\newcommand{\Supp}[0]{{\operatorname{Supp}}}
\newcommand{\Bs}[0]{{\operatorname{Bs}}}

\newtheorem{thm}{Theorem}[section]
\newtheorem{lem}[thm]{Lemma}

\newtheorem{prop}[thm]{Proposition}

\theoremstyle{definition}
\newtheorem{defn}[thm]{Definition}
\newtheorem{rem}[thm]{Remark}
\newtheorem{ack}{Acknowledgments}        
\newtheorem*{notation}{Notation}         %\renewcommand{\thenotation}{} 

\theoremstyle{remark}

\begin{document}
\bibliographystyle{amsalpha+}
\maketitle

\abstract 
Let $(X,\Delta)$ be a proper dlt pair 
and $L$ a nef Cartier divisor such that $aL-(K_X+\Delta)$ is nef and 
log big on $(X,\Delta)$ for some $a\in {\mathbb Z}_{>0}$. 
Then $|mL|$ is base point free for every $m\gg 0$. 
\endabstract

\setcounter{section}{-1}
\section{Introduction}\label{se1}
The purpose of this paper is to prove the following theorem. 
This type of base point freeness was suggested by M.~Reid in \cite 
[10.4]{R}.  

\begin{thm}[Base point free theorem of Reid-Fukuda type]
\label{th} 
Assume that $(X,\Delta)$ is a proper dlt pair. 
Let $L$ be a nef Cartier divisor such that $aL-(K_X+\Delta)$ is nef and 
log big on $(X,\Delta)$ for some $a\in {\mathbb Z}_{>0}$. 
Then $|mL|$ is base point free for every $m\gg 0$, 
that is, there exists a positive integer $m_0$ such that 
$|mL|$ is base point free for every $m\geq m_0$. 
\end{thm}

This theorem was proved by S.~Fukuda in the case where 
$X$ is smooth and $\Delta$ is a reduced simple normal crossing divisor 
in \cite {fk2}. 
In \cite {fk3}, he proved it on the assumption that 
$\dim X \leq 3$ by using log Minimal Model Program. 
Our proof is similar to \cite {fk3}. 
However, we do not use log Minimal Model Program even in $\dim X \leq 3$. 
He also treated this problem under some extra conditions in \cite {fk4}. 

\begin{ack}
I would like to express my sincere gratitude to 
Dr.~Daisuke Matsushita for giving me some comments. 
\end{ack}

\begin{notation}
(1) We will make use of the standard notations and definitions 
as in \cite{KM}.   

(2) A pair $(X,\Delta)$ denotes that $X$ is a normal variety over $\mathbb C$ 
and $\Delta$ is a $\bQ$-divisor on $X$ such that $K_X+\Delta$ is 
$\bQ$-Cartier. 

(3) $\xDiff$ denotes the different (See \cite [Chapter 16] {FA}).   
\end{notation}

\section{Preliminaries}

In this section, we make some definitions and collect the necessary 
results. 

\begin{defn}(cf. \cite [Definition 1.3]{Ka2})
A subvariety $W$ of $X$ is said to be a {\em {center of log canonical 
singularities}} for the pair 
$(X,\Delta)$, if there exists a proper birational morphism 
from a normal variety $\mu: Y\to X$ and a prime divisor $E$ on $Y$ 
with the discrepancy $a(E,X,\Delta)\leq -1$ such that $\mu(E)=W$.  
\end{defn}

\begin{defn}
Let $(X,\Delta)$ be lc and $D$ a $\bQ$-Cartier $\bQ$-divisor on $X$. 
The divisor $D$ is called {\em {nef and log big}} on $(X,\Delta)$ 
if $D$ is nef and big, and $(D^{\dim W}\cdot W)>0$ 
for every center of log canonical 
singularities $W$ of the pair 
$(X,\Delta)$.
\end{defn}

\begin{rem}
(1) Our definition of nef and log big is equivalent to 
that of Reid and Fukuda (See \cite[Definition]{fk3}). 

(2) In \cite{fj}, the center of log canonical singularities 
of a dlt pair was investigated 
(See \cite [Definition 4.6, Lemma 4.7]{fj}). 
\end{rem}

The following proposition is \cite [Proposition 2]{fk3} 
(for the proof, see \cite [Proof of Theorem 3]{fk2} and \cite 
[Lemma 3]{Ka}). 

\begin{prop}\label{prop}
Let $(X,\Delta)$ be a proper dlt pair and 
$L$ a nef Cartier divisor such that $aL-(K_X+\Delta)$ is nef 
and big for some $a\in {\mathbb Z}_{>0}$. 
If $\Bs |mL| \cap \llcorner \Delta\lrcorner=\emptyset$ for 
every $m\gg 0$, then $|mL|$ is base point free for every $m\gg 0$, 
where $\Bs |mL|$ is the base locus of $|mL|$. 
\end{prop}

\begin{lem}\label{vanishing}{\em{(cf. \cite [Lemma]{fk1})}}
Let $X$ be a proper smooth variety and $\Delta=\sum_{i} d_i \Delta_i$ 
a sum of distinct prime divisors 
such that $\Supp \Delta$ is a simple normal crossing divisor and 
$d_i$ is a rational number with $0\leq d_i\leq 1$ for every $i$. 
Let $D$ be a Cartier divisor on $X$. 
Assume that $D-(K_X+\Delta)$ is nef and log big on $(X,\Delta)$. 
Then $H^{i} (X ,{\mathcal O} _{X} (D))=0$ for every $i>0$. 
\end{lem}

This is a generalization 
of Kawamata-Viehweg vanishing theorem.

\section{Proof of Theorem}
\proof[Proof of Theorem (\ref {th})]
By using \cite[Resolution Lemma]{S} as in the proof of the 
Divisorial Log Terminal Theorem of \cite{S}, we have a log resolution 
$f:Y\to X$ of $(X,\Delta)$, which satisfies the following conditions; 
\begin{enumerate}
\item [(1)] $K_Y+f^{-1}_{*}\Delta=f^{*}(K_X+\Delta) +\sum _{i} a_i E_i$ 
with $a_i>-1$ for every $i$, where $E_i$'s are irreducible exceptional 
divisors,  
\item [(2)] $f$ induces isomorphism at every generic point of center 
of log canonical singularities of $(X,\Delta)$. 
\end{enumerate}
We define $E:=\sum _{i}\ulcorner a_i \urcorner E_i\geq 0$ and 
$F:=f^{-1}_{*}\Delta+E-\sum _{i}a_i E_i$. 
Then $K_Y+F=f^{*}(K_X+\Delta) +E$. 
If $\llcorner\Delta\lrcorner=0$, then $(X,\Delta)$ is klt. 
So we may assume that $\llcorner\Delta\lrcorner\neq 0$. 
We take an irreducible component $S$ of 
$\llcorner\Delta\lrcorner$. 
Then $(S,\xDiff (\Delta-S))$ is dlt. 
It can be checked easily by \cite [Corollary 5.52, 
Definition 2.37]{KM} and \cite[17.2 Theorem]{FA}. 
We put $S_0:=f^{-1}_{*}S$ and $M:=f^{*}L$. 
We consider the following exact sequence; 
$$
0\to {\mathcal O}_{Y} (-S_0) \to {\mathcal O}_{Y}
\to {\mathcal O}_{S_0}\to 0.
$$ 
Tensoring with ${\mathcal O}_{Y}(mM+E)$ for $m\geq a$, 
we have the exact sequence; 
$$
0\to {\mathcal O}_{Y} (mM+E-S_0) \to {\mathcal O}_{Y}(mM+E)
\to {\mathcal O}_{S_0}(mM+E)\to 0.
$$ 
By Lemma (\ref{vanishing}), $H^{1}(Y,{\mathcal O}_{Y} (mM+E-S_0))=0$. 
Note that $M$ is nef and 
$mM+E-S_0-(K_Y+F-S_0) = f^{*} (mL-(K_X+\Delta))$ 
is nef and log big on $(Y,F-S_0)$. 
Then we have that 
$$
H^{0}(Y,{\mathcal O}_{Y} (mM+E))\to H^{0}(S_0, {\mathcal O}_{S_0} 
(mM+E))
$$ 
is surjective. By the projection formula, we have that 
$$
H^{0}(Y,{\mathcal O}_{Y} (mM+E))\simeq H^{0}(X,f_{*}{\mathcal O}_{Y} (mM+E)) 
\simeq H^{0}(X,{\mathcal O}_{X} (mL)) 
$$ 
and  
$$
H^{0}(S_0, {\mathcal O}_{S_0}(mM+E))\supset 
H^{0}(S_0, {\mathcal O}_{S_0}(mM))\simeq 
H^{0}(S,{\mathcal O}_{S}(mL)). 
$$
Note that $E$ is effective and $f$-exceptional and that $E|_{S_0}$ 
is effective but not necessarily $f|_{S_0}$-exceptional, 
where $f|_{S_0}: S_0\to S$. 
We consider the following commutative diagram; 
$$
\begin{CD}
H^{0}(Y,{\mathcal O}_{Y}(mM+E))
@>>> H^{0}(S_0,{\mathcal O}_{S_0}(mM+E))@>>>0\\
@AA\text{$\cong$}A @AA\text{$\iota$}A\\
H^{0}(X,{\mathcal O}_{X}(mL)) @>>> H^{0}(S,{\mathcal O}_{S}(mL)).
\end{CD}
$$
Then $H^{0}(X,{\mathcal O}_{X}(mL))\to H^{0}(S,{\mathcal O}_{S}(mL))$ is 
surjective and $\iota$ is isomorphism 
since the left vertical arrow is isomorphism and 
$\iota$ is injective by the above argument. 
By induction on dimension, $|mL|_{S}|$ is base point free 
for every $m\gg 0$ since $(aL-(K_X+\Delta))|_S = 
aL|_{S} -(K_S+\xDiff (\Delta -S))$ is nef and log big on 
$(S,\xDiff (\Delta-S))$. 
So we have that $\Bs |mL|\cap \llcorner\Delta\lrcorner=\emptyset$. 
By Proposition (\ref {prop}), we get the result. 
\endproof 

\ifx\undefined\bysame
\newcommand{\bysame}{\leavevmode\hbox to3em{\hrulefill}\,}
\fi

\end{document}